\documentclass[12pt,reqno]{amsart}
\usepackage[latin1]{inputenc}  
\usepackage{amssymb}
\usepackage{stmaryrd}
\usepackage[all]{xy}
\usepackage{graphicx}
\usepackage{hyperref}
\theoremstyle{plain}

\newtheorem{theore}{Theorem}
 
\newtheorem{theo}{Theorem}[section]
\newtheorem{lemme}[theo]%
 {Lemma}
{Definition}
\newtheorem{prop}[theo]%
{Proposition}
\newtheorem{coro}[theo]%
{Corollary}
{Definition-Proposition}
\newtheorem{conjecture}[theo]%
{Conjecture}

\newcommand{\finpreuve}{\mbox{} \hfill \mbox{$\Box$}}
\newcommand{\finpreuvelemme}{\mbox{} \hfill \mbox{$\blacksquare$}}
\newenvironment{preuve}{\noindent {\it Proof :}}{\finpreuve}
\newenvironment{preuvelemma}{\noindent {\it Proof of the lemma:
}}{\finpreuvelemme}

\makeatletter
\newcounter{enum@ux}
\def\brkenum#1{%
\setcounter{enum@ux}{\value{enum\romannumeral\the\@enumdepth}}%
\end{enumerate} #1%
\begin{enumerate}%
\addtocounter{enum\romannumeral\the\@enumdepth}{\value{enum@ux}}}
\makeatother

\makeindex

\newcommand{\pg}{\mathfrak{p}}
\newcommand{\Pg}{\mathfrak{P}}
\newcommand{\qg}{\mathfrak{q}}

\newcommand{\Q}{\mathbb{Q}}
\newcommand{\FF}{\mathbb{F}_r(t)}
\newcommand{\KK}{\mathcal{K}}
\newcommand{\rg}{\mathfrak{r}}

\title{On Tsfasman--Vl\u adu\c t Invariants of Infinite Global Fields}
\author{Philippe Lebacque}
\date{}
\address{NWF I - Mathematik, Universität Regensburg, D-93040 Regensburg
Germany}
\email{Philippe.Lebacque@mathematik.uni-regensburg.de}
\begin{document}
\maketitle

\begin{abstract} In this article we study certain asymptotic properties of global fields. We consider the set  of Tsfasman--Vl\u adu\c t invariants of infinite global fields and answer some natural questions arising from their work. In particular, we prove the existence of infinite global fields having finitely many strictly positive invariants at given places, and the existence of infinite number fields with certain prescribed invariants being zero. We also give precisions on the deficiency of infinite global fields and on the primes decomposition in those fields.
\end{abstract}
\section{Introduction}

\renewcommand{\thetheo}{\Alph{theo}}

In the 80', Ihara \cite{IPD} initiated the asymptotical theory of global fields in the particular case of unramified infinite global fields. More recently, Tsfasman and Vl\u adu\c t  \cite{TVF}  generalised his work to any infinite global fields, and defined a set of invariants of such fields.  The aim of this paper is to investigate their properties further.

\subsection{} Following \cite{TVF}, let us first recall some definitions and some basic facts about the theory of infinite global fields. Let $r$ be a power of a prime number $p.$ An infinite global field $\mathcal{K}$ is an infinite separable algebraic extension of $\Q$ or $\FF$ such that $\mathcal{K}\cap\bar{\mathbb{F}_r}=\mathbb{F}_r.$ In the function field case, all the extensions we consider are separable and without constant extensions. We write $(NF)$ (respectively $(FF))$ in order to signify that an assertion concerns the number field case (resp. the function field case). For any global field $K,$ put $n_K$ its degree over $\Q$ or $\FF,$ $g_K$ its genus ($\frac{1}{2}\log{|D_K|}$ in the number field case). Denote by $P(K)$ (resp. $P_f(K))$ its set of places (resp. of non-archimedean places). For any place let $\pg\in P_f(K),$ $\mathrm{N}\pg$ be its norm. For any prime power $q,$ let $\Phi_q(K)$ denote the number of its places of norm $q,$ and let $\Phi_\mathbb{R}(K),$ $\Phi_\mathbb{C}(K)$ be the number of its real and complex places respectively. Let $A$ denote the set of parameters $\{\mathbb{R},\ \mathbb{C},\ p^k, \ p \text{ prime number, } k\in\mathbb{N}^\ast\}$ in the number fields case, $\{ r^k,\ k\in \mathbb{N}^\ast\}$ in the function fields case, and let $A_f$ be the subset of the parameters which are prime powers. We also define relative $\Phi$-numbers as follows: given a place $\pg$ of a global field $E,$ and $K/E$ a finite extension, let $\Phi_{\pg,q}(K)$ denote the number of places of $K$ above $\pg$ with the norm $q$ (or which are real or complex in the case $q=\mathbb{R}$ or $\mathbb{C}$).  We will omit the field $K$ in our notation if there is no possible confusion.

Given a set of primes $S$ of a global field $K,$ denote by $\delta(S)$ its Dirichlet density if it exists:
$$\delta(S)=\lim_{s\to 1^{+}}\frac{\sum_{\pg\in S}\mathrm{N}\pg^{-s}}{\sum_{\pg\in P_f(K)}\mathrm{N}\pg^{-s}}.$$ If it does not exist, we put $$\bar{\delta}(S)=\limsup_{s\to 1^{+}}\frac{\sum_{\pg\in S}\mathrm{N}\pg^{-s}}{\sum_{\pg\in P_f(K)}\mathrm{N}\pg^{-s}}.$$ For any set of places $S$ of a global field $K$ and any extension of global fields $L/K$ (resp. $K/L$), let $S(L)$ be the set of places of $L$ lying above (resp. under) $S.$

 Recall that a sequence $\{K_i\}_{i\in\mathbb{N}}$ of global fields is said to be a family if $K_i$ is not isomorphic to $K_j$ for $i\neq j,$ and if  the constant field of all the $K_i$ equals one and the same $\mathbb{F}_r$ (see \cite{TVF}). In any family of global fields, we have $g_i\to\infty,$ because for any real number $g_0,$ there is only finitely many number fields (respectively function fields  over a given constant field up to isomorphism) whose genus does not exceed $g_0.$

 We say that a family of global fields $\{K_i\}_{i\in\mathbb{N}}$ is asymptotically exact if the limit $$\phi_q(\{K_i\})=\lim_{i\to\infty} \frac{\Phi_q(K_i)}{g_{K_i}}$$ exists for any $q\in A.$ A sequence $\{K_i\}_{i\in\mathbb{N}}$ of global fields is called a tower if $K_i$ is strictly included in $K_{i+1}$ for any integer $i.$  Any tower of global fields $\{K_i\}$ is an asymptotically exact family, and the $\phi_q$'s depend only on $\cup_{i\in\mathbb{N}} K_i.$ Therefore one can define Tsfasman-Vl\u adu\c t invariants $\phi_q(\mathcal{K})$ of an infinite global field $\KK$ as being the $\phi_q$'s corresponding to any tower $\{K_i\}_{i\in\mathbb{N}}$ such that $\displaystyle\cup_{i\in\mathbb{N}} K_i=\mathcal{K}.$ Define the support of an infinite global field as $$\mathrm{Supp}(\KK)=\{q\in A,\ \phi_q(\KK)> 0\}.$$ We also define the quantity $$\phi_\infty(\KK)=\lim_{i\to\infty}\frac{n_{K_i}}{g_{K_i}}.$$

We say that an infinite global field is asymptotically good if its support is not empty, and asymptotically bad if it is empty. In the same manner we can define the relative invariants $\phi_{\pg,q}(\KK),$ for a global field $E,$ a place $\pg$ of $E$ and an infinite global field $\KK$ containing $E,$ as the limit of the ratio $\Phi_{\pg,q}(K_i)/g_{K_i}$ for any tower $\{K_i\}$ such that $E\subset K_i$ for any $i$ and $\cup_{i\in\mathbb{N}} K_i=\KK.$ As there are only finitely many places of $E$ whose norm is not greater than $q$ and the number of archimedean places is also finite,  we have $$\phi_q(\mathcal{K})=\sum_{\pg\in P(E)}\phi_{\pg,q}(\mathcal{K}).$$ Define also the prime support of $\KK/E$ as the set $$\mathrm{PSupp}(\mathcal{K}/E)=\{\pg\in P_f(E)\ | \ \exists\, q\in A_f \ \phi_{\pg,q}>0 \}.$$
 
 Finally, for any (possibly infinite) global field $\KK,$ any global field $K\subset\KK,$ let $Ram(\KK/K)$ (respectively $Rams(\KK/K)$, resp. $Dec(\KK/K)$) denote the ramification locus (resp. the wild ramification locus, resp. the decomposition locus) of $\KK/K,$ that is set of places of $K$ that are ramified in some (resp. that are wildly ramified in some, resp. that are split in any) finite subextension $L/K$ of $\KK/K.$  In the case $K=\Q$ or $\FF,$ we omit it in the notation.

 \subsection{}Tsfasman and Vl\u adu\c t proved that an infinite number field  is asymptotically good if and only $\phi_\infty>0.$ In the function field case that is necessary but not sufficient.
 More precisely, for any prime number $p,$  $$(NF)\quad \sum_{m\geq 1}m\phi_{p^m}\leq\phi_\mathbb{R}+2\phi_\mathbb{C}.$$ There is an analogous property for function fields considering instead of the $\Phi_{r^m}$'s the numbers $\Phi_{\pg,r^m}$'s of places above a given place $\pg$ of $\FF.$
 
Except for the case of asymptotically bad infinite global fields or of those constructed from the optimal ones (see \cite{GSO}), we do not know any example where the set of the invariants of an infinite global field or even its support is completely known. But Tsfasman and Vl\u adu\c t also gave the following inequalities for the invariants which generalize Drinfeld-Vl\u adu\c t bound (see \cite{TVF}, Theorems ):
 
 \begin{theore}[Tsfasman-Vl\u adu\c t Basic Inequalities]\label{ineg}  For any infinite global field, we have the following inequalities:
 \begin{align*}
(NF-GRH)\text{  }& \sum_q\frac{\phi_q\log{q}}{\sqrt{q}-1} + (\log{\sqrt{8\pi}}+\frac{\pi}{4}+\gamma/2)\phi_\mathbb{R}+(\log{8\pi}+\gamma)\phi_\mathbb{C}\leq 1,\\
(NF)\text{  }\ & \sum_q\frac{\phi_q\log{q}}{q-1}+(\gamma/2+\log{2\sqrt{\pi}})\phi_\mathbb{R}+(\gamma+\log{2\pi})\phi_\mathbb{C}\leq 1,\\
(FF)\text{  }\ & \sum_{m=1}^\infty \frac{m\phi_{r^m}}{r^{\frac{m}{2}}-1}\leq 1,
\end{align*}
where the (GRH) indication means once and for all that an assertion is true assuming the Generalised Riemann Hypothesis.
 \end{theore}

Several questions arise naturally. For any family $(a_q)_{q\in A}$ of real numbers satisfying the basic  inequalities, does it exist an infinite global field such that $\phi_q=a_q$ for all $q\in A?$ We will see that the answer to this question is negative in the number field case, and the analogous result holds for function fields (see Corollary \ref{kunyavski}). Some other natural weaker questions are far from being within reach at the moment. For example, we do not know if  there exists an infinite global field with an infinity of $\phi_q$ being positive, neither if there are infinite number fields with all but one invariants equal to zero (the function field case being known in the case where $r$ is a square, using optimal towers). Even if we are not able to give answers to this two questions, we can prove the following result:
\begin{theo}\label{phi1}
Let $n$ be an integer and $t_1,...,t_n\in A_f.$ There exists an infinite global field (both in the number field and function field cases) such that $\phi_{t_1},...,\phi_{t_n}$ are all $>0,$ and such that, in the number field case, any other $\phi_q$, with $\mathrm{gcd}(q,\prod t_i)>1$ is zero.
\end{theo}
Note that  an analogous result holds for function fields, if we consider the $\Phi_{\pg,r^m}$ numbers instead of $\Phi_{r^m}.$

\begin{theo}\label{phi2} Let $P$ be a finite set of prime numbers. There exists an asymptotically good infinite Galois number field $\mathcal{K}$ such that, for any positive integer $m$ and any $p\in P,$ $\phi_{p^m}(\mathcal{K})=0.$ Moreover, $\mathrm{PSupp}(\KK/\Q)$ equals $Dec(\KK/\Q)$ (and therefore has a zero Dirichlet density, see Proposition \ref{zero}).
\end{theo}

\subsection{}Tsfasman and Vl\u adu\c t also defined the deficiency of an infinite global field as the difference between the two terms of the basic inequalities:
\begin{align*}
(NF-GRH)\text{  }\   \delta^{(1)}=1 &-\sum_q\frac{\phi_q\log{q}}{\sqrt{q}-1} 
\\ & - (\log{\sqrt{8\pi}}+\frac{\pi}{4}+\gamma/2)\phi_\mathbb{R}-(\log{8\pi}+\gamma)\phi_\mathbb{C},\\
(NF)\text{  }\  \delta^{(2)}=1 &-\sum_q\frac{\phi_q\log{q}}{q-1}\\
& -(\log{2\sqrt{\pi}}+\gamma/2)\phi_\mathbb{R}-(\log{2\pi}+\gamma)\phi_\mathbb{C},\\
(FF)\text{  }\  \delta^{(3)}=1 &-\sum_{m=1}^\infty \frac{m\phi_{r^m}}{r^{\frac{m}{2}}-1}.
\end{align*}

In the case of function fields, when $r$ is a square, there are infinite function fields which reach the deficiency zero, and they are called optimal (see \cite{GS} for example). But in the case of number fields, such fields are unknown and we think that their existence is doubtful. 

The deficiency is an increasing function, therefore optimal fields, if they exist, have to satisfy some similarly properties to those satisfied by just-infinite global fields (meaning that they have no proper infinite subextension):
\begin{theo} \label{defaut} For any $i\in \{1,2,3\}$, the map $\mathcal{K}\mapsto\delta^{(i)}(\mathcal{K})$ is an increasing map for the inclusion of infinite global fields. 
\end{theo}

\subsection{} One can refine in some sense Theorem \ref{phi1}, looking for the best possible field having $n$ distinct invariants positive, in term of deficiency.  
\begin{theo} There exists an infinite number field having $n$ positive non archimedean invariants, such that its deficiency $\delta_n\leq 1-\varepsilon_n$ satisfies:
\begin{align*}
(NF-GRH) \quad \varepsilon_n\sim &\frac{8}{3\sqrt{n\log{n}}} \text{ and,}\\
(NF) \quad \varepsilon_n\sim& \frac{4}{3n}.
\end{align*}
\end{theo}

One can obtain a corresponding result for function fields, but the class field theory would  likely give a very bad estimation comparing to what we could obtain starting with an optimal tower. This should be done in a further work.
\subsection{} As one can see from the definition of the ${\phi_q}'s,$ these invariants are closely related to the decomposition of primes in infinite global fields. The most general result concerning it is a corollary of the Cebotarev density theorem. 
\begin{prop}\label{zero} Let $\KK$ be an infinite number field (resp. function field), and let $T$ be the set of places of $\mathbb{Q}$ (resp. $\FF)$ that split in $\KK.$ Then $\delta(T)$ exists and is equal to $0.$
\end{prop}
This result implies in particular that, for an infinite Galois global field $\KK$ over a global field $K,$ the set $$\{\pg\in P(K)\  | \ \phi_{\pg,\mathrm{N}\pg}>0\}$$ has a zero Dirichlet density, and in the particular case of Galois number fields, the set $\{p\ |\ p \text{ is a prime number and } \phi_p>0\}$ has to be very small.
It seems to be hard to improve it, but in his paper \cite{IPD}, Ihara considered the primes decomposition in infinite unramified Galois extensions of global fields, and proved that, for any number field $k$ and any unramified Galois extension $\KK$ of $k$  
$$\sum_{\pg\in P_f(k)} \frac{\log{\mathrm{N}\pg}}{\mathrm{N}\pg^{f(\pg)}-1}\leq C(k),$$ where $f(\pg)$ is the inertia index of a place of $\KK$ above $\pg,$ and $C(k)$ is a real number depending only on $k.$ 
Using the Tsfasman-Vl\u adu\c t-basic inequalities, one deduces the following:
\begin{prop} \label{iharaz} Let $\KK$ be an infinite number field (resp. function field) and let $T$ denote the set of places of $\mathbb{Q}$ (resp. $\FF)$ that split in $\KK.$ 
Then \begin{align*}(NF) \quad &\sum_{\mathfrak{p}\in T}\frac{\log{\mathrm{N}(\mathfrak{p})}}{\mathrm{N}(\mathfrak{p})-1}\leq \frac{1}{\phi_\infty},\\
(FF) \quad &\sum_{\mathfrak{p}\in T}\frac{\log_r{\mathrm{N}(\mathfrak{p})}}{\sqrt{\mathrm{N}(\mathfrak{p})}-1}\leq \frac{1}{\phi_\infty}, 
\end{align*} where $1/\phi_\infty$ can be $+\infty$ (corresponding to asymptotically bad cases).
\end{prop}

However in the case of asymptotically bad infinite global field, the sum $s(T)=\sum_{\pg\in T}\mathrm{N}\pg^{-1}$ can be infinite, even if the set of ramification is relatively small:

\begin{theo}\label{ex} There exists an infinite global field $\KK$ (in both cases of number and function fields) without wild ramification such that $s(Dec(\KK))=\infty$ and $\delta(Ram(\KK))=0.$ \end{theo}
One can improve this result a bit, giving a more precise information on the ramification locus: for any $\varepsilon>0$ there is an infinite global field $\KK$ satisfying the above properties, and such that $s(Ram(\KK))\leq\varepsilon.$

\medskip

The structure of the paper is as follows. First we recall basic facts concerning infinite global fields, which will be useful for the comprehension of the remainder parts. In $\S$ $3$ we prove Theorem \ref{phi1} using class field towers. The following paragraph is devoted to the study of the deficiency and the proof of Theorem \ref{defaut}. In $\S$ $5,$ we estimate the defect of towers involved in Theorem \ref{phi1}. After that, we consider the problem of prime decomposition in infinite global fields, which is central in our study and prove Theorem \ref{ex}. Finally, we prove Theorem \ref{phi2} in the last section.

The author would like to thank Michael Tsfasman for all he has done as his advisor,  Christian Maire for the time he spent answering his questions and for his advice concerning \S \ref{schmi}, and Alexei Zykin for all of his comments on the first version of this paper.

\setcounter{theo}{0}
\renewcommand{\thetheo}{\arabic{section}.\arabic{theo}}

\section{Basic facts on infinite global fields}

In this paragraph we recall briefly some basic properties of asymptotically good infinite global fields. For more details see \cite{Lth}, or \cite{GS} for the the function field Galois case. If we want to construct an asymptotically good infinite global field, we ask for three conditions. First, the tower should be tamely ramified, or we should be able to control the wild ramification which is often the most difficult part. Second, the tower should be unramified outside of a finite set of primes. There are examples in the function field case, where $\phi_\infty>0$ and where the set of ramification is infinite (in fact everywhere ramified, see \cite{PO}), but we do not know any example where the field is asymptotically good. And third, there should be at least one split prime (in the function fields case). 

Controlling the wild ramification leads to specific calculations, as you can find in \cite{GSO}. In the following, we will always avoid wild ramification. If there is no deep wild ramification, and if the field is Galois, then an asymptotically good has to be finitely ramified.
We say that an infinite global field $\KK$ is almost-Galois if there is a tower $\{K_i\}_{i\in\mathbb{N}}$ of global fields such that $\cup K_i=\KK$ and $K_{i+1}/K_i$ is Galois for any $i\in\mathbb{N}.$ 

 \begin{prop}\label{fini}
Let $\mathcal{K}$ be an infinite number field (resp. function field) such that $$\inf_{K' \text{finite}}\#Rams(\mathcal{K}/K')=0.$$ Then $\mathcal{K}$ is asymptotically good (resp. $\phi_\infty(\mathcal{K})>0$) if its ramification locus $S$ is finite. Moreover, if $\mathcal{K}/K$ (for a global field $K$) is Galois, the converse is true. If $\mathcal{K}=\cup K_i$ is almost-Galois (for a tower $\{K_i\}_{i\in\mathbb{N}}$ of global fields), then $\mathcal{K}$ is asymptotically good (resp. $\phi_\infty(\mathcal{K})>0$) if and only if $$\sum_{i=0}^\infty \frac{1}{n_{K_i}}\sum_{\mathfrak{p}\in Ram(K_{i+1}/K_i)}\log\mathrm{N} \mathfrak{p}<\infty,$$ where $\log$ denote the base $e$ (resp. base $r$) logarithm.
\end{prop}
\begin{preuve} This is a straightforward application of the Riemann-Hurwitz formula (see \cite{Lth}). For a number field (resp. a function field) $K,$ let $g_K^\ast$ denote $g_K$ (resp. $g_K-1$), so that the Hurwitz formula is exactly the same in the case of number fields and of function fields without constants extensions.

Let $\{K_i\}_{i\in\mathbb{N}}$ $(K_0=K)$ be a tower representing $\mathcal{K}$ (meaning that $\cup K_i=\KK)$. Suppose that the ramification locus of $\KK$ is finite. Let $K_j$ be a global field of genus $>1$ such that $\#Rams(\mathcal{K}/K_j)=0$. We can suppose that $j=0$ without any loss of generality. We have to show that the ratio $g_{K_i}^\ast/n_{K_i}$ is bounded. We can always replace $n_{K_i}$ by $[K_i:K_0],$ the computations being modified by a factor equal to the degree of $K_0$ over the ground field ($\Q$ or $\FF).$

 Write $g_i^\ast, n_i$ for $g_{K_i}^\ast, [K_i:K_0]$ respectively. Let $S$ be $Ram(\mathcal{K}/K_0)$ and $S_i=Ram(K_i/K_0)$. The Riemann-Hurwitz formula (see \cite[1.3.10]{XNB} and \cite[III,3.13]{NEU}) implies that $$ g_i^\ast=n_i g_0^\ast+\frac{1}{2}\sum_{\mathfrak{p}\in S_i}\sum_{\mathfrak{P}/\mathfrak{p}} (e_\mathfrak{P}-1)f_{\mathfrak{P}}\log \mathrm{N} \mathfrak{p},$$
where $\mathfrak{P}$ is a place above $\mathfrak{p}$, $e_\mathfrak{P}$ denote the ramification degree of $\mathfrak{P}/\mathfrak{p}$ and $f_\mathfrak{P}$  the inertia degree of $\mathfrak{P}/\mathfrak{p}$.

As our extensions are separable, $$ \sum_{\mathfrak{P}/\mathfrak{p}} e_\mathfrak{P}f_{\mathfrak{P}}=n_i, $$ and one deduces

\begin{align*}
g_i^\ast & \leq n_i g_0^\ast+\frac{1}{2} n_i\sum_{\mathfrak{p}\in S_i}\log\mathrm{N} \mathfrak{p},\\
& \leq n_i g_0^\ast+\frac{1}{2} n_i\sum_{\mathfrak{p}\in S}\log\mathrm{N} \mathfrak{p},
\end{align*}
this last sum being finite because $S$ is finite, and the first statement is proven. 

Suppose now that $S$ is infinite and that the tower $\{K_i\}$ is Galois, so that $K_i/K_0$ is Galois. We have to show that $g_i^\ast/n_i$ is not bounded. 

$$ g_{i}^\ast=[K_{i},K_0] g_0^\ast+\frac{1}{2}\sum_{\mathfrak{p}\in S_i}\sum_{\mathfrak{P}/\mathfrak{p}} (e_\mathfrak{P}-1)f_{\mathfrak{P}}\log \mathrm{N} \mathfrak{p}.$$
As $e_\mathfrak{P}-1\geq e_\mathfrak{P}/2$ because $e_\mathfrak{P}\geq 2$ (the extension is Galois), we deduce that

\begin{align*}
g_{i}^\ast & \geq n_i g_0^\ast+\frac{1}{2}\sum_{\mathfrak{p}\in S_i}\log \mathrm{N} \mathfrak{p}\sum_{\mathfrak{P}/\mathfrak{p}} \frac{e_\mathfrak{P}}{2}f_{\mathfrak{P}}\\
&\geq n_i g_0^\ast+\frac{1}{4}n_i \sum_{\mathfrak{p}\in S_i}\log \mathrm{N} \mathfrak{p}.
\end{align*}
The sum is now divergent because $\#S_i$ tends to $+\infty$, which prove the converse.

Let us prove now the assertion concerning the almost-Galois case. Let $S_i$ denote the ramification locus of the Galois extension $K_{i+1}/K_i$. Applying the Riemann-H\"urwitz formula as before, we get:
\begin{align*}
\frac{n_{i+1}}{n_i} g_i^\ast+\frac{1}{4}\frac{n_{i+1}}{n_i}\sum_{\mathfrak{p}\in S_i}\log\mathrm{N} \mathfrak{p} & \leq g_{i+1}^\ast \ \leq \frac{n_{i+1}}{n_i} g_i^\ast+\frac{1}{2}\frac{n_{i+1}}{n_i}\sum_{\mathfrak{p}\in S_i}\log\mathrm{N} \mathfrak{p}\\
\frac{g_i^\ast}{n_i}+\frac{1}{4\,n_i}\sum_{\mathfrak{p}\in S_i}\log\mathrm{N} \mathfrak{p} & \leq\frac{g_{i+1}^\ast}{n_{i+1}} \leq \frac{g_i^\ast}{n_i}+\frac{1}{2\,n_i}\sum_{\mathfrak{p}\in S_i}\log\mathrm{N} \mathfrak{p}.
\end{align*}
We obtain by induction $$\frac{g_0^\ast}{n_0}+\frac{1}{4}\sum_{j=0}^i\sum_{\mathfrak{p}\in S_j}\frac{\log\mathrm{N} \mathfrak{p}}{n_j} \leq \frac{g_{i+1}^\ast}{n_{i+1}} \leq \frac{g_0^\ast}{n_0}+\frac{1}{2}\sum_{j=0}^i\sum_{\mathfrak{p}\in S_j}\frac{\log\mathrm{N} \mathfrak{p}}{n_j},$$ which concludes the proof.
\end{preuve}

In particular, any infinite global field $\KK$ containing a global field $K$ such that $\KK/K$ is unramified,  satisfies $\phi_\infty>0.$ The following result explains why the split places are of particular interest in the study of the $\phi_q$'s.

\begin{prop} \label{propdec} Let $\mathcal{K}$ be an infinite global field, and let $K$ be a global field contained in $\KK.$
 Suppose that $\phi_\infty>0$ and that a non-archimedean place $\pg$ of $K$ splits in $\KK/K.$ Then $\phi_{\mathrm{N}\pg}>0.$ Moreover, if $\mathcal{K}$ is Galois, and if there is a $q\in A_f$ such that $\phi_q>0,$ then $\phi_\infty>0$ and there is a global field $L\subset\KK$ containing $K$ and a non-archimedean place $\pg$ of $K$ such that any place of $L$ above $\pg$ is of norm $q$ and splits in $\KK/L.$
 \end{prop}

\begin{preuve} For the sake of notation, let us prove it in the number fields case. If $\pg$ is split in $\KK/K$ then, for any tower $\{K_i\}$ representing $\KK$ such that $K_0=K,$ $\Phi_{\mathrm{N}{\pg}}(K_i)\geq n_{K_i}/[K_0:\Q],$ which implies the first assertion. Suppose now that $\KK/K$ is Galois, and that $\phi_q>0,$ for a prime power $q.$  There is a place $\pg$ of $K$ such that $\phi_{\pg,q}>0,$ because $\sum_{\pg\in P(K)}\phi_{\pg,q}=\phi_q.$ Then $\phi_\infty\geq\phi_{\pg,q}$ has to be non zero. All the places above $\pg$ in $\KK$ have the same ramification index and inertia degree, which have to be both finite: indeed, if $K'$ is a global field contained in $\KK,$ the ramification index $e,$ the inertia degree $f,$ and the number $d$ of places above $\pg$ satisfy $e\,f=[K':K]/d,$ and this ratio is bounded over the subfields $K'$ of $\KK$ because of our assumption. Therefore all the places above $\pg$ have to split in $\KK/L$ for some global field $L.$
\end{preuve}

\section{On class field towers}

\subsection{}This section is devoted to the construction of infinite global fields using unramified classfield towers. For definition and results concerning $S$-ramified $T$-split class field towers, see \cite{MA1}. Let us recall first the construction: given a global field $K,$ a prime number $p$ and a set $T$ of places of $K,$ let $H_\ell^T(K)$ denote the maximal unramified abelian $\ell$-extension of $K$ where $T$ is split. Construct a tower of field as follows: let $K_0$ be our global field $K$ and $T_0$ a set of non archimedean places of $K_0.$ For $i\geq 1$ put $K_i=H_\ell^{T_{i-1}}(K_{i-1}),$ and let $T_i$ be the set of places of $K_i$ above $T_{i-1}.$ The tower $\{K_i\}$ is called the $\ell$-$T_0$-class field tower of $K_0.$ The question of the finiteness of this tower remained open for years until Golod and Schafarevitch gave  in the 1960's a criterion to prove that such a tower is infinite. This result in particular implies the following theorem:
\begin{theo}\label{TVNX}[Tsfasman--Vl\u adu\c t \cite{TVF} (NF), Serre \cite{SLN} , Niederreiter--Xing \cite{XNB} (FF)]
Let $K/k$ be a cyclic extension of global fields of degree $\ell.$ Let $T(k)$ be a finite set of non archimedean places of $k$ and $T(K)$ be the set of places above $T(k)$ in $K.$ Suppose in the function field case that $\mathrm{gcd}\{\ell,\deg\pg,\pg\in T(K)\}=1.$ Let $Q$ be the ramification locus of $K/k.$ Let 
\begin{align*} 
(FF)\quad C(T,K/k)=& \#T(k)+2+\delta_\ell+2\sqrt{\#T(K)+\delta_\ell},\\
(NF)\quad C(T,K/k)=& \#T(K)-t_0+r_1+r_2+\delta_\ell+2-\rho+\\
& 2\sqrt{\#T(K)+\ell(r_1+r_2-\rho/2)+\delta_\ell},
\end{align*}
where $\delta_\ell=1$ if $K$ contains the $\ell$-root of unity, and $0$ otherwise, $t_0$ is the number of principal ideals in $T(k),$ $r_1=\Phi_\mathbb{R}(K),$ $r_2=\Phi_\mathbb{C}(K)$ and $\rho$ is the number of real places of $k$ which become complex in $K.$
Suppose that 
$\#Q\geq C(T,K/k).$ Then $K$ admits an infinite unramfied $\ell$-$T(K)$-class field tower.
\end{theo}
Remark that the assumption on the degree of the places in $T$ in the function field case guaranties that there is no constant field extension in the tower.
These unramified towers give us infinite global fields such that $\phi_{\mathrm{N}\pg}>0$ for every $\pg\in T,$ provided we can construct sufficiently ramified cyclic extensions. Even though this point can be made explicitly (see \S\ref{explicit}), we will use the Grunwald--Wang theorem.

\subsection{} Let $S$ be a set of primes of a global field $k,$ containing archimedean places in the number field case. Let $k_S$ denote the maximal extension of $k$ unramified outside of $S,$ which is Galois over $k.$ Put $G_S=Gal(k_S/k).$ We set $$(NF)\quad \quad\mathbb{N}(S)=\{n\in \mathbb{N}\ |\ v_\mathfrak{p}(n)=0\text{ for any }\mathfrak{p}\notin S\}$$ in the number fields case, and let $\mathbb{N}(S)$ denote the set of numbers prime to the characteristic of $k$ in the function fields case.

Let us recall that the exponent of a finite group $\#A$ is the smallest integer $a$ such that $x^a=1_A$ for every $x$ in $A.$
 
Let us now state the Grunwald--Wang theorem, as we can find it in \cite{NCG} in the case where $S$ contains all the places of $k.$

\begin{theo}[Grunwald--Wang] Let $S$ be a set of places of a global field $k,$ containing the archimedean primes in the number field case, such that $\delta(S)=1,$ let $T\subset S$ be a finite set of primes of $k$ and $A$ be a finite abelian group of exponent $a$ such that $\#A\in\mathbb{N}(S)$. Let $K_\mathfrak{p}/k_\mathfrak{p}$ be, for any $\mathfrak{p}\in T,$  local abelian extensions such that $G(K_\mathfrak{p}/k_\mathfrak{p})$ can be imbedded in $A$. Then there is a global abelian extension $K/k$ with Galois group $A,$ unramified outside of $S$ such that $K$ has the given completions $K_\mathfrak{p}$ for any $\mathfrak{p}\in T$, with the exception of the special case ($k,a,T$) when the following four conditions hold:
\begin{enumerate}
\item $a=2^ra', \ \ r> 2$
\item $k$ is a number field
\item $k(\mu_{2^r})/k$ is not cyclic
\item $\left\{\mathfrak{p}\in S\ \  |\ \mathfrak{p} \text{ divides } 2\text{ and } \mathfrak{p} \text{ is not split in } k(\mu_{2^r})/k\right\}\subset T.$
\end{enumerate}
\end{theo}

\begin{preuve} It is nothing but the proof of $9.2.3.$ of \cite{NCG}, where $k$ is replaced by $k_S.$ For the sake of completeness, let us recall it. Consider $A$ as a trivial $G_S$-module. First recall the following.
For any $\mathfrak{p}\in S,$ we choose a $k$-embedding $k_S\to \bar{k}_\mathfrak{p}$ by chosing a place $\bar{\mathfrak{p}}$ of $k_S$ above $\mathfrak{p}.$ This induces the restriction map from $G(\bar{k}_\mathfrak{p}/k_\mathfrak{p})$ to $G_S=G(k_S/k)$ (whose image is the decomposition group of  $\bar{\mathfrak{p}}$ over $k$). It induces the map $$H^1(G_S,A)\to H^1( G(\bar{k}_\mathfrak{p}/k_\mathfrak{p}),A)$$ independent from the choice of the embedding. 

$$\xymatrix@1{
H^1(k_S|k,A)\ar[r]^<<<<{res}  & \prod_{\pg\in T} H^1(\bar{k}_\mathfrak{p}|k_\mathfrak{p},A).
}$$

Let us prove the theorem now. We want to show that the map $$\mathrm{Epi}(G_S,A)\to \prod_T\mathrm{Hom}(G(\bar{k}_\mathfrak{p}/k_\mathfrak{p}),A)$$ is onto, where $\mathrm{Epi}(G_S,A)$ denotes the onto morphisms from $G_S$ to $A.$ Let $\mathfrak{q}_1,...,\mathfrak{q}_r$ be finite places of $S-T$ not dividing $2$ in the number field case, and let $$\psi_{\mathfrak{q}_i}: G(\bar{k}_{\mathfrak{q}_i}/k_{\mathfrak{q}_i})\to A$$ be morphisms such that their images generate $A.$ For any $\mathfrak{p}\in T,$ $\psi_\mathfrak{p}:G(\bar{k}_{\mathfrak{p}}/k_{\mathfrak{p}})\to A$ denotes the canonical map induced by the chosen embedding $G(K_\mathfrak{p}/k_\mathfrak{p})\to A.$ Let $T'$ denote the set $T\cup\{\mathfrak{q}_1,...,\mathfrak{q}_r\}.$   
The map 
$$\xymatrix{
H^1(k_S|k,A)\ar@{->}[r]^<<<<{res}&\prod_{T'} H^1(\bar{k}_\mathfrak{p}|k_\mathfrak{p},A)
}$$
is known to be onto in the case where $\delta(S)=1$ and $A$ is a trivial $G_S(k)$-module, with exception of the special case $(k,exp(A),T')=(k,exp(A),T)$ (see \cite[9.2.2]{NCG}).
 An inverse image  $\psi:G_S\to A$ of $(\psi_{\mathfrak{p}})_{\mathfrak{p}\in T'}$ in $H^1(G_S,A)=\text{Hom}(G_S,A)$ realises the local extensions and is onto because of the choice of the $\psi_{\mathfrak{q}_i}.$
\end{preuve}

\subsection{}Using this theorem, we can now prove the following important result:

\begin{coro}\label{tourbonne} Let $P$ be a finite set of places of $\Q$ (resp. of $\mathbb{F}_r(t)$, containing at least one rational place). Then there is an asymptotically good infinite number field $\KK$ (resp. infinite function field) such that all the places in $P$ are split, and that there exists a global field $K$ of prime degree $2$ such that $\KK/K$ is unramified.
\end{coro}
\begin{preuve} Let $Q$ be a finite set of places of non archimedean places of $\mathbb{Q}$ (resp. of $\mathbb{F}_r(t)$) such that $P\cap Q=\emptyset,$  and $\#Q$ is big enough to satisfy the condition of theorem \ref{TVNX} for any quadratic extension of $\Q$ (resp. $\FF$). 

  Applying the theorem with $A=\mathbb{Z}/2\mathbb{Z}$, $T=P\cup Q,$ we obtain a cyclic extension of degree $2$ of $\mathbb{Q}$ (resp. de $\mathbb{F}_r(t)$) where the set $Q,$ at least, is ramified, and such that $P$ is split. Note that we can obtain the totally ramified local extensions of degree $2$ by adjoining a root of Eisenstein separable polynomials. Then we use Theorem \ref{TVNX} to obtain an unramified $2$-class field tower where all the places in $P$ are split. This tower does not have constants extension because $P$ contains a rational place, and is asymptotically good because of the results of $\S 2.$
\end{preuve}

We remark that we can choose a tower with $\phi_\mathbb{R}$ or $\phi_\mathbb{C}$ positive, choosing in the first step a real quadratic field or an imaginary one. In the corollary, we can also take towers of degree $p$ different from the characteristic of our fields, but in that case, the norm of the places constituting $Q$ must equal $1$ modulo $p$ so that the totally ramified local $p$-extensions exist.

\subsection{Proof of Theorem \ref{phi1}}


In the function field case, we apply the corollary \ref{tourbonne} with a set of places $P$ containing at least one place of degree $\log_r{t_i}$ for any $i=1,...,n,$ which gives directly the desired result. If we want to obtain such additional properties as in the number field case, we just follow the number field case proof. 
 
In the number field case, we consider the set $P=\{p^{(1)},\dots,p^{(k)}\}$ of prime numbers $p^{(1)},\dots,p^{(k)}$ which divide one of the $t_i.$

\begin{lemme} There is a finite extension $L$ of $\mathbb{Q}$ having the following properties:
\begin{enumerate}\item for any $i=1\dots n$ there exists $\pg\in P(L)$ such that $\mathrm{N}\pg=t_i,$
\item for any $\pg\in P(L),$ if there exists $p\in P$  such that $p|\mathrm{N}\pg$ then there is $i\in\{1,...,n\}$ such that $\mathrm{N}\pg=t_i.$
\end{enumerate}
\end{lemme}

Assuming this lemma, let us apply the corollary \ref{tourbonne} to the set of places $P.$ It gives us an infinite number field $\KK=\cup K_i$ represented by the tower of number field $\{K_i\},$ $K_0=\mathbb{Q},$ unramified over $K_1,$ where all the places of $P$ are split. Then the tower $\{L.K_i\}_{i\in\mathbb{N}}$ is unramified over $L.K_1,$ therefore asymptotically good, and all the places of $L$ above $P$ are split in $L.K_i.$ So we obtain the field $L.\KK$ having the desired properties. 

\begin{preuvelemma}
We construct the tower $L_0=\mathbb{Q}\subset L_1,...\subset L_k=L$ using the Grunwald--Wang theorem. Let us begin by ordering the set $T$.  Write $$T=\left\{ t_1^{(1)},...,t_{i_1}^{(1)};t_1^{(2)},...,t_{i_2}^{(2)};t_1^{(k)},...,t_{i_k}^{(k)} \right\}$$ such that $p^{(r)}$ divide all  $t_1^{(r)},...,t_{i_r}^{(r)}$ but no others.
Put $d_j^{(r)}=\log_{p^{(r)}}t_j^{(r)}.$ 
 
 Recall first the following properties of local fields guarantees that we can use the Grunwald--Wang theorem (see \cite[$\S$ III$.5,$ Theorem $2$]{CL}):
\begin{lemme} Let $K$ be a local field with a finite residue field $k,$ and let $k'/k$ be a cyclic extension of $k.$ Then there exists an unramified cyclic extension $K'/K$ such that $k'$ is the residue field of $K'.$ 
\end{lemme}
Given a global field $K$ and a place $\pg,$ the completion $K_\pg$ of $K$ at $\pg,$ an integer $n>0,$
there exists an unramified cyclic extension $K'$ of $K_\pg,$ of degree $n$ (meaning that $\pg$ is totally inert in it).

Using the Grunwald--Wang theorem, let us construct by induction a tower $L_0\subset\dots\subset L_k$ having the following properties: 
\begin{enumerate}\item in $L_{r+1}/L_r,$ all the places above $p^{(j)}$ are split for $j\neq r+1,$
\item there exists an extension $L_r^1$ of $L_r$ of degree $i_{r+1},$ contained in $L_{r+1}$ such that all the places above $p^{(r+1)}$ are split.
\item in $L_{r+1}/L_r^1$ there are $[L_{r+1}:L_0]/(i_{r+1}d^{(r+1)}_j)$ places of norm $t_j^{(r+1)}$ for any $j.$
\end{enumerate}

$L_0$ is given. Let us suppose that $L_r/L_0$ has been constructed. Put $m_r=[L_r:L_0].$ Let $L_r^1/L_r$ be an extension of  degree $i_{r+1}$ where all the places above any $p^{(i)}$ are split. In $L_r^1$ one has got $i_{r+1}m_r$ places $$p_{1,1}^{(r+1)},\dots,p_{1,m_r}^{(r+1)},\dots,p_{i_{r+1},1}^{(r+1)},\dots,p_{i_{r+1},m_r}^{(r+1)}$$ above $p^{(r+1)},$ put into $i_{r+1}$ packets containing $m_r$ places each.

Consider then successive cyclic extensions  of prime degree (so we do not fall in the special case). We obtain an extension $L_{r+1}/L_r^1$ such that:
\begin{enumerate}
\item all the places above $p^{(i)}$ are split for any $i\neq r+1,$
 \item above $p_{j,l}^{(r+1)},$ for any $j$ and any $l,$ there are $[L_{r+1}/L_r^1]/d^{(1)}_j$ places of norm $t^{(r+1)}_j.$ 
 \end{enumerate}
In order to deal with $j^\text{th}$ packet, we ask for the existence of a cyclic extension of prime degree dividing $t_j^{(r+1)},$ in which all the places above the $p_{j,l}^{(r+1)}$ are totally inert for any $l$, and in which all the other pointed places are split, until we obtain places of norm $t_j^{(r+1)}.$

By induction, we obtain therefore $L_k=L,$ satisfying by the construction the two given properties. This ends the proof of the lemma and that of Theorem \ref{phi1}.
 \end{preuvelemma}
 
\section{Deficiency and optimal fields}\label{defi}

The deficiency of an infinite global field is defined as the difference between the two sides of the basic inequality. Tsfasman  and Vl\u adu\c t proved that it is related to the limit distribution of the zeroes of the zeta function. This distribution admits a continuous density, and the deficiency is its value at $0$ (see \cite{TVF}). The study of the deficiency is therefore not only simpler than the study of the invariants itself, but also gives us some very interesting knowledge about the field.

\subsection{Proof of Theorem \ref{defaut}}

The first but not the least result concerning the deficiency is Theorem \ref{defaut}. Let us prove it in the case of $\delta^{(1)},$ the proof being exactly the same in the two other cases (we replace $\log$ by $\log_r$ in the function field case). We begin by treating the non-archimedean terms:
\begin{lemme} For any prime number $p,$ any $m\in\mathbb{N},$ and any infinite global fields $\mathcal{K}\subset\mathcal{L}$, $$(NF)\qquad \sum_{k=1}^m\frac{k\phi_{p^k}(\mathcal{L})\log{p}}{\sqrt{p^k}-1}\leq\sum_{k=1}^m\frac{k\phi_{p^k}(\mathcal{K})\log{p}}{\sqrt{p^k}-1},$$
$$(FF)\qquad \sum_{k=1}^m\frac{k\phi_{r^k}(\mathcal{L})}{\sqrt{r^k}-1}\leq\sum_{k=1}^m\frac{k\phi_{r^k}(\mathcal{K})}{\sqrt{r^k}-1}.$$
\end{lemme}
Remark that we could replace $\log{p}/(p^{k/2}-1)$ by any decreasing function. We can also prove the same inequality for the $\phi_{\pg,q}$ along the same lines.

\begin{preuvelemma}
Let $\mathcal{K}\subset\mathcal{L}$ be two infinite number fields and $p$ a prime number. Recall that, for any $m,$ we have: 
$$A_m(\mathcal{K}):= \sum_{k=1}^m k\phi_{p^k}(\mathcal{K})\geq A_m(\mathcal{L}).$$ Using Abel transform on $\sum_{k=1}^m\frac{k\phi_{p^k}(\mathcal{K})\log{p}}{\sqrt{p^k}-1},$ we get: 
\begin{align*}
\sum_{k=1}^m\frac{k\phi_{p^k}(\mathcal{K})\log{p}}{\sqrt{p^k}-1}= & A_{m}(\mathcal{K})\frac{\log{p}}{\sqrt{p^{m}}-1}\\
& +\sum_{k=1}^{m-1}A_k(\mathcal{K})\left(\frac{\log{p}}{\sqrt{p^{k}}-1}-\frac{\log{p}}{\sqrt{p^{k+1}}-1}\right),
\end{align*}
from which we deduce the lemma. \end{preuvelemma}

Taking the limit, we obtain, for two infinite global fields $\mathcal{K}\subset\mathcal{L}:$
$$(NF)\quad\sum_{k=1}^\infty\frac{k\phi_{p^k}(\mathcal{L})\log{p}}{\sqrt{p^k}-1}\leq\sum_{k=1}^\infty\frac{k\phi_{p^k}(\mathcal{K})\log{p}}{\sqrt{p^k}-1},$$
$$(FF)\qquad \sum_{k=1}^\infty\frac{k\phi_{r^k}(\mathcal{L})}{\sqrt{r^k}-1}\leq\sum_{k=1}^\infty\frac{k\phi_{r^k}(\mathcal{K})}{\sqrt{r^k}-1}.$$
Let us now consider the archimedean terms:
\begin{lemme} \label{lemmearchi} For two infinite number fields $\mathcal{K}\subset\mathcal{L}$ and two positive real numbers $\alpha_1,$ $\alpha_2$ such that $2\alpha_1\geq\alpha_2,$ one has: $$\alpha_1 \phi_\mathbb{R}(\mathcal{L})+\alpha_2 \phi_\mathbb{C}(\mathcal{L}) \leq \alpha_1\phi_\mathbb{R}(\mathcal{K})+\alpha_2\phi_\mathbb{C}(\mathcal{K}).$$
\end{lemme}
\begin{preuvelemma} Let us consider two towers of number fields $\mathcal{K}=\cup K_i,$ $\mathcal{L}=\cup L_i,$ with $K_i\subset L_i,$ and let us write the relations between archimedean places of $K_i$  and $L_i.$ 

We forget the indices and denote by $K$ and $L$ these fields. Let $P_\mathbb{R}(K)$ denote the set of real places of $K.$ Recall that the complex places are always split in $L/K$, giving birth to $[L:K]$ complex places of $L$ and each real place $v$ of $K$ decomposes into $\Phi_{v,\mathbb{R}}(L/K)$ real places and $\Phi_{v,\mathbb{C}}(L/K)$ complex places of $L,$ so that $\Phi_{v,\mathbb{R}}(L/K)+2\Phi_{v,\mathbb{C}}(L/K)=[L:K].$ 

Therefore one has $$\Phi_{\mathbb{R}}(L)=\sum_{v\in P_\mathbb{R}(K)} \Phi_{v,\mathbb{R}}(L/K),$$ and $$\Phi_{\mathbb{C}}(L)=[L:K]\Phi_{\mathbb{C}}(K)+\sum_{v\in P_\mathbb{R}(K)}\Phi_{v,\mathbb{C}}(L/K).$$
Let $\alpha_1$ and $\alpha_2$ be real numbers such that $2\alpha_1\geq\alpha_2.$
Then
\begin{align*}\alpha_1 \Phi_{\mathbb{R}}(L)+\alpha_2 \Phi_{\mathbb{C}}(L)&=\alpha_1\sum_{v\in P_\mathbb{R}(K)}\Phi_{v,\mathbb{R}}(L/K)+\alpha_2[L:K]\Phi_{\mathbb{C}}(K)\\ 
&+\alpha_2\sum_{v\in P_\mathbb{R}(K)}\Phi_{v,\mathbb{C}}(L/K)\\
&\leq \alpha_1\sum_{v\in P_\mathbb{R}(K)}\left(\Phi_{v,\mathbb{R}}(L/K)+2\,\Phi_{v,\mathbb{C}}(L/K)\right)\\
\\ & +\alpha_2[L:K]\Phi_{\mathbb{C}}(K)\\
&\leq [L:K]\big(\alpha_1\Phi_{\mathbb{R}}(K)+\alpha_2\Phi_{\mathbb{C}}(K)\big).
\end{align*}
As $g_L\geq [L:K]g_K,$ one obtains, for $g_K>0$:
$$\alpha_1 \frac{\Phi_{\mathbb{R}}(L)}{g_L}+\alpha_2 \frac{\Phi_{\mathbb{C}}(L)}{g_L}\leq \alpha_1\frac{\Phi_{\mathbb{R}}(K)}{g_K}+\alpha_2\frac{\Phi_{\mathbb{C}}(K)}{g_K}.$$ Taking the limit we get the lemma.\end{preuvelemma}

This two lemmas (the first in the function field case) give us directly, taking for $\alpha_1$ and $\alpha_2$ corresponding to the archimedean terms in the deficiency, the decreasing property of the map  $\mathcal{K}\mapsto 1-\delta^{(i)},$ for any $i,$ which ends the proof.

\subsection{Optimal Fields}

We say that an infinite global field is optimal for $\delta^{(i)}$ if its deficiency $\delta^{(i)}$ is equal to $0.$ In the function case, when $r$ is a square, there are examples of optimal fields. Different constructions can be used, such as tower of modular curves or recursive towers. In the number field case, or even in the function field case where $r$ is not a square, the question whether optimal towers exist or not remains open. Let us give first some properties that should be satisfied by optimal fields. We will prove then that most of the infinite number fields constructed using the class field theory cannot be optimal.

\begin{prop} \label{gal1} Let $\mathcal{K}$ be an optimal infinite number field (resp. function field) for  $\delta^{(i)}.$ If $\pg$ is a non archimedean place of $\mathbb{Q}$ (resp. $\FF$) such that $\phi_{\pg,q}=0$ for any $q,$ there is no infinite global field $\mathcal{L}$ contained in $\KK$ such that $\phi_{\pg,q}(\mathcal{L})>0.$
\end{prop}

\begin{preuve} Let us prove it for $i=1.$ Suppose there exists such an extension $\mathcal{L}.$ The proof of Theorem \ref{defaut} shows that, for any prime number $p,$ $$\delta_p(\mathcal{L}):=\sum_{k=1}^\infty\frac{k\phi_{p^k}(\mathcal{L})\log{p}}{\sqrt{p^k}-1}\geq\sum_{k=1}^\infty\frac{k\phi_{p^k}(\mathcal{K})\log{p}}{\sqrt{p^k}-1}=\delta_p(\mathcal{K}).$$ Because of Lemma \ref{lemmearchi}, one has
$$1-\big(\sum_{p\neq\ell}\delta_{p}(\mathcal{L})+\delta_\ell(\mathcal{L})+\alpha_1 \phi_\mathbb{R}(\mathcal{L})+\alpha_2 \phi_\mathbb{C}(\mathcal{L})\big)\leq 1-(1+\delta_\ell(\mathcal{L})),$$ where $\alpha_1$ and $\alpha_2$ are the archimedean coefficients involved in the deficiency $\delta^{(1)}.$
 Therefore $\delta^{(1)}(\mathcal{L})\leq -r\phi_{\ell^r}\log\ell/(\ell^{r/2}-1)<0,$ which leads to a contradiction. \end{preuve}

\begin{prop} \label{gal2} Let $\mathcal{K}$ be an optimal infinite Galois number field (resp. function field) for $\delta^{(i)}.$ Let $\pg$ be a place of $\Q$ (resp. $\FF$). Suppose there is a prime number $p$ and an integer $d$ such that $\phi_{\pg,p^d}(\KK)>0.$  Then $\mathcal{K}$ does not contain any infinite global subfield $\mathcal{L}$ such that $\phi_{\pg,p^m}(\mathcal{L})>0,$ $m\neq d.$ 
\end{prop}

\begin{preuve} We show the result for $\delta^{(1)}.$ Let $\mathcal{L}\subset \mathcal{K}$ such that one $\phi_{\pg,p^m}>0,$ and $m<d$ (the other case cannot happen because $\mathcal{K}$ is Galois).  We prove that $\mathcal{L}$ is a subfield of $\mathcal{K}$ whose deficiency is strictly smaller than the deficiency of $\mathcal{K}.$ Indeed, for any prime $q\neq p,$ we have $\delta_q(\mathcal{L})\geq\delta_q(\mathcal{K}).$
As for $\pg,$ one has
$$\sum_{k\leq d}k\phi_{p^k}(\mathcal{L})\geq d\phi_{p^d}(\mathcal{K}),$$ since $\mathcal{L}\subset \mathcal{K},$ we deduce $$\sum_{k\leq d}\frac{k\phi_{p^k}(\mathcal{L})}{\sqrt{p^k}-1}>\sum_{k\leq d}\frac{k\phi_{p^k}(\mathcal{L})}{\sqrt{p^d}-1}\geq \frac{d\phi_{p^d}(\mathcal{K})}{\sqrt{p^d}-1}.$$
But this is impossible, because $\KK$ is optimal.
\end{preuve}

Let us prove now that fields constructed with class field theory are mainly not optimal. Consider the deficiency without $GRH$ in the number field case (the $GRH$-deficiency result can easily be deduced), and call it $\delta.$

Let $\mathcal{Q}$ denote $\Q$ in the number field case (resp. $\FF$ in the function field case). Let $K/\mathcal{Q}$ be a finite Galois extension. Let $S(K)$ be a finite set of non-archimedean places of $K.$ Let $P$ be a minimal set of non-archimedean places of $K$ such that:
\begin{enumerate} \item $P$ is disjoint from $S(K),$
\item $P$ is stable under the action of the Galois group of $K/\mathcal{Q},$
\item For any non-archimedean places $\pg$ and $\qg$ of $K$ outside of $S(K),$ such that $\mathrm{N}\pg<\mathrm{N}\qg,$ if $\qg$ belongs to $P,$ $\pg$ belongs to $P.$
\item In the number fields case, $$(NF)\quad\sum_{\mathfrak{p}\in P}\frac{\log{\mathrm{N} \mathfrak{p}}}{\mathrm{N} \mathfrak{p}-1}>g^\ast(K)-\frac{\ell\alpha_2}{2},$$ where $\alpha_2=\gamma+\log{2\pi},$ and in the function fields case, $$(FF)\quad\sum_{\mathfrak{p}\in P}\frac{\log_r{\mathrm{N} \mathfrak{p}}}{\sqrt{\mathrm{N} \mathfrak{p}}-1}>g^\ast(K).$$
\end{enumerate}
 
The sum $\sum_{\mathfrak{p}\in P}\frac{\log{\mathrm{N} \mathfrak{p}}}{\mathrm{N} \mathfrak{p}-1}$ (resp. $\sum_{\mathfrak{p}\in P}\frac{\log_r{\mathrm{N} \mathfrak{p}}}{\sqrt{\mathrm{N} \mathfrak{p}}-1}$) taken over all non archimedean places of $K,$ is divergent; therefore such $P$ exists. One constructs it taking consecutively all the non archimedean places of $K$ outside of $S(K)$ (and those obtained by applying the Galois action), until the sum becomes bigger than the right hand side term. Let  $\pg_0$ be a prime of maximal norm in $P,$ and $p_0$ the prime of $\mathcal{Q}$ under it. Put
  \begin{align*}
  (NF)\quad \alpha(K,S):=&\frac{\ell\log p_0}{g^\ast(K)}\left(\frac{1}{\mathrm{N} \pg_0-1}-\frac{1}{\mathrm{N} \pg_0^\ell-1}\right)\\
 (FF)\quad \alpha(K,S):=&\frac{\ell\log_r\mathrm{N} p_0}{g^\ast(K)}\left(\frac{1}{\sqrt{\mathrm{N} \pg_0}-1}-\frac{1}{\sqrt{\mathrm{N} \pg_0^\ell}-1}\right)
 \end{align*}

Let us now state the theorem.

\begin{theo} \label{nonopt} Let $K/\mathcal{Q}$ be a cyclic extension of prime degree $\ell$ ramified exactly at $S$ a finite set of non archimedean places of $\mathcal{Q}.$ Let $\qg,\qg'\notin S$ be two places of $K$ (of  relatively prime degrees in the function field case).  Let $\KK={K}_\emptyset(\ell)$ be the maximal unramified $\ell$-extension of $K$ (resp. let $\KK=K_\emptyset^{\{\qg,\qg'\}}(\ell)$ be the maximal unramified $\ell$-extension of $K$ where $\qg$ and $\qg'$ are split). Suppose that $\#S\geq 3+\ell+2\sqrt{2+\ell^2}$ (resp.  $\#S\geq 6-\varepsilon+2\sqrt{3\,\ell+1-\varepsilon} $--- see Theorem \ref{TVNX} for notation).

Then $\KK$ is not optimal. Moreover its deficiency satisfies $\delta(\KK)\geq \alpha(K,S).$ 
\end{theo}
Remark that it is also the case of any field containing $\KK,$ in particular it is the case of $\Q_S$ (the maximal extension of $\Q$ unramified outside of $S$) in the number field case. The condition on $S$ guarantees that $\KK$ is infinite because of  Theorem \ref{TVNX}.

\begin{preuve} Suppose $\delta(\KK)<\alpha.$ Let $\log$ denote as usual the logarithm  with base $e$ in number field case, and with base $r$ in the function field case. Prove first that there is a finite place with no contribution in $\delta.$
\begin{lemme} There is a place $\pg_1\in P$ such that, for any $m>0,$ $\phi_{p_1^m}=0,$ where $p_1=\pg_1\cap\Q.$ 
\end{lemme}
\begin{preuvelemma}
Suppose that for every place $\pg$ in $P,$ there is an $m_\pg>0$ such that $\phi_{\pg,\mathrm{N} \mathfrak{p}^{m_\pg}}>0.$ As $\KK$ is Galois, it is also the case for all the places above $\pg\cap\mathcal{Q}.$ $S$ is by our hypothesis sufficiently large so that, for any place $\rg$ not in $S,$ the maximal extension $\KK^\rg$ of $K,$ contained in $\KK,$ with all the places above $\rg$ in $K$ split, is not finite. 
Let us prove first that, for any place $\pg$ of $P,$ $m_\pg=1.$ If we only want to prove the non optimality, it is sufficient to apply the proposition \ref{gal2}.

Let us suppose that $m_\rg>1$ for a place $\rg.$ All the places above $\rg$ in $K$ have the same norm, they are unramified in any tower contained in $\KK,$ therefore the difference of their contributions to the deficiency of $\KK^\rg$ and their contributions to that of $\KK$ is given by 
\begin{align*}
(NF)\quad  &\frac{\ell\log \mathrm{N}\rg}{g^\ast(K)({\mathrm{N} \mathfrak{p}}-1)}-\frac{\ell\log  \mathrm{N}\rg}{g^\ast(K)({\mathrm{N} \mathfrak{r}}^{m_\rg}-1)}.\\
(FF)\quad &\frac{\ell\log_r  \mathrm{N}\rg}{g^\ast(K)({\sqrt{\mathrm{N}} \mathfrak{p}}-1)}-\frac{\ell\log  \mathrm{N}\rg}{g^\ast(K)(\sqrt{{\mathrm{N} \mathfrak{r}}^{m_\rg}}-1)}.
\end{align*}
It cannot exceed $\alpha,$ because the deficiency belongs to $[0,1].$ Moreover, this quantity is decreasing in $p,$ increasing in $m,$ so it is sufficient to verify that, for the biggest $p$ in our set, and for the smallest $m$ (meaning $m=\ell$), the condition is satisfied.

Because of the definition of $\alpha,$ this is not satisfied and therefore we have a contradition.
Consequently $m_\pg=1$ for any place $\pg$ of $P.$

Then for any $\pg\in P$, $$\phi_{\mathrm{N} \pg}=\frac{\Phi_{\mathrm{N} \pg}(K)}{g^\ast(K)}.$$ Indeed, all the places above $\pg\cap\Q$ in $K$ are split and the equality follows.
 
We have 
\begin{align*}(NF)\quad \sum_q\phi_q\frac{\log{q}}{q-1}+\delta_\infty \geq&\frac{1}{g^\ast(K)}\sum_{\pg\in P}\frac{\log{\mathrm{N}\pg}}{\mathrm{N} \pg-1}+\frac{\ell\alpha_2}{2\,g^\ast(K)}>1,\\ 
(FF)\quad \sum_q\phi_q\frac{\log_r{q}}{\sqrt{q}-1} \geq&\frac{1}{g^\ast(K)}\sum_{\pg\in P}\frac{\log_r{\mathrm{N}\pg}}{\sqrt{\mathrm{N} \pg}-1}>1,
\end{align*}
where $\delta_\infty=\alpha_1 \phi_\mathbb{R}(\mathcal{L})+\alpha_2 \phi_\mathbb{C}(\mathcal{L})$ is the contribution to the deficiency of the archimedean places.
 Indeed $\delta_\infty\geq\frac{1}{2}\alpha_2\phi_\infty$ and $\phi_\infty=\ell\,g^\ast(K)^{-1}.$ This contradicts T-V Basic Inequalities.
 \end{preuvelemma}

Consider now the maximal unramified extension of $K$ such that $\pg_1$ is split. The contribution of $\rg_1=\pg_1\cap \mathcal{Q}$ to the deficiency of this extension is $\frac{\ell\log p_1}{g^\ast(K)(\mathrm{N} \pg_1-1)}>\alpha(K,S),$ whereas it is zero in $\KK.$ This leads to a contradiction since this infinite global fields would have a strictly negative deficiency.
\end{preuve}

\section{An Effective Version of Theorem \ref{phi1}}\label{explicit}

The aim of this paragraph is to produce an example of infinite number field with $n$ positive invariants having a deficiency $\delta_n$ as small as possible. 

In order to achieve our goal, we will apply Theorem \ref{TVNX}. Take for the set P the first prime $n$ numbers greater than 2: $P=\{3,5,7,...,p_n\}$
We will take for $K$ a quadratic extension, $\mathbb{Q}(\sqrt{q_1...q_{r_n}r})$, where $r_n$ is the smallest integer satisfying the conditions of Theorem \ref{TVNX} with $s=2n:$  
$$r_n=1+\lfloor{n+5+2\sqrt{2n+4+1}}\rfloor.$$
Then $r_n=n+2\sqrt{2n} +c_n$ where $c_n$ is a bounded sequence. Remark that for imaginary quadratic extensions, we could take a smaller $r_n,$ but that does not change anything asymptotically.

Now let us choose $r.$ We want the $p_i$ to be split (otherwise they are inert and we loose a  $\sqrt{n}$ factor). Take for $q_i$, $i \leq r_n$, the $r_n$ primes following $p_n$ in the prime progression, and choose $r$ in the following way so that the $p_i$ are split:

\begin{lemme}\label{temppp}
There exists $r\in \{0,...,2\prod p_i \prod q_j-1\}$ prime to $2, p_i,q_j$ for any $i$, $j$ such that the $p_i$ are split in $K/\mathbb{Q}$. 
\end{lemme}

\begin{preuvelemma} Recall that the $p_i$ are split in $K/\mathbb{Q}$ if $r\prod q_j $ is a square mod $p_i$ for any $i,$ i.e. if for any $i\leq n$ $$\left(\frac{r\prod q_j}{p_i}\right)=\left(\frac{r}{p_i}\right)\prod_j \left(\frac{q_j}{p_i}\right)= 1.$$
When the $p_i$ and the $q_j$, $j\leq r_n$ have been chosen,  $r$ has to satisfy
$$\forall i\in\{1,...,n\} \ \ \ \ \left(\frac{r}{p_i}\right)=\prod_{j}\left(\frac{q_j}{p_i}\right).$$

We choose a non zero solution modulo $p_i$ for any $i.$ We lift it to an odd integer  of $\{0,...,2\prod p_i \prod q_j-1\}$ equal to $1$ modulo $q_j$ for any $j.$ If this integer is different from $1,$ it splits into a product of prime numbers, which are distinct from the $p_i$, $q_i$ and from $2.$ We choose it squarefree, which is possible because if $r$ has a square factor $p^2,$ $r/p^2$ suits us too.
Note that $r$ can be equal to $1$, and that its factors ramify in $K/\mathbb{Q}$ but these factor have got no consequences on the validity of Theorem \ref{TVNX}.
\end{preuvelemma}
We remark that we could probably improve the upper bound on $r.$ But it would only improve the involved constants in the following proposition. For instance we took $r$ in the most trivial way. We could have computed the number of integers smaller than $2.\prod p_i$ satisfying the conditions on the Legendre's symbols and compared it to the number of integers prime to all the $q_i$ and $2.$

From now on we denote by $K_n$ the field $K$ corresponding to $n$ and an $r$ obtained from Lemma \ref{temppp}. 
\begin{prop}\label{genre}
The genus $g_{K_n}$ of $K_n$ satisfies 
$$g_{K_n}\leq g_n, \text{ where } \ \ g_n=\frac{3}{2}n \log{n} +o(n\log{n}),$$
$$g_{K_n}\geq g'_n, \text{ where } \ \ g'_n=\frac{1}{2}n \log{n} +o(n\log{n}).$$
\end{prop}

\begin{preuve}
$$g_{K_n}=\frac{1}{2}\sum_j \log{q_j}+\frac{1}{2}\log{r}+\frac{1}{2}\varepsilon\log{2},$$
where $\varepsilon=2$ if $2$ is ramified and $0$ otherwise.
One has $\log{r}\leq \sum_i \log{p_i}+\sum_j \log{q_j}+\log{2}.$ Therefore

 $$\log{r}\leq \sum_{p\leq q_{r_n}}\log p.$$

Thus $$g_{K_n} \leq \sum_{p\leq q_{r_n}}\log{p} - \frac{1}{2}\sum_{p\leq p_{n}}\log{p}+\log 2=g_n.$$

From the analytic number theory (see \cite{HW}), we know  that the function $\vartheta(x)=\sum_{p\leq x} \log p$ is asymptotically equal to $x+o(x),$ as $x\to\infty,$ and that the $n^{\text{th}}$ odd prime number $p_n$ satisfies $p_n=n\log{n}+o(n\log{n})$. 
Thus $$\sum_{p\leq p_{n}}\log{p}=n\log{n}+o(n\log{n}).$$

As $r_n=n+2\sqrt{2n}+\mathcal{O}(1),$ one has $q_{r_n}=(2n+2\sqrt{2n})\log{n}+o(n\log{n})$, $q_{r_n}$ being the $(n+r_n)^\text{th}$ odd prime number.
Therefore: 
$$\sum_{p\leq q_{r_n}}\log{p}=(2n+2\sqrt{2n})\log{n}+o(n\log{n}),$$
and we get $$g_n=\frac{3}{2}n \log{n} +o(n\log{n}).$$
The second inequality can be obtained in the same way, using $r\geq 1.$ 
\end{preuve}

Apply now Theorem \ref{TVNX} to $K_n/\mathbb{Q}$. $K_n$ admits an infinite unramified class field tower $\KK_n=(K_{n}^i)_{i\geq 0} ,$ with $K_{n}^0=K_n$, where the $p_i$ are split in $\KK_n/\Q$ for $i\leq n.$
 
\begin {prop} For any $i\geq 0,$ 
$g_{K_n^i}=2^i g_{K_n}$
\end{prop}
\begin{preuve} 
As the tower is unramified, we have for any $i$, $$g_{K_{n}^{i+1}}=g_{K_{n}^i}[K_{n}^{i+1}:K_{n}^i].$$ We obtain the desired result by induction.
\end{preuve}
  
\begin{coro}
$\phi_\infty(\KK_n)=O((n\log{n})^{-1}).$
\end{coro}
\begin{preuve}
For any $i$, we have $n_{K_n^i}=2^{i+1}$, so that
$$\frac{n_{K_n^i}}{g_{K_n^i}}=\frac{2^{i+1}}{g_{K_n^i}}=\frac{2}{g_{K_n}}$$
and, taking the limit we get $$\phi_\infty(\KK_n)=\frac{2}{g_{K_n}}.$$
Therefore $\phi_\infty(\KK_n)\leq \frac{2}{g_n'}$ and the result follows from Proposition \ref{genre}.
\end{preuve}

Let us evaluate now the sum involved in the deficiency.  Because the $p_i$ are split and because of Proposition \ref{genre}, we have  $\phi_{p_i}=2\frac{1}{g_{K_n}}.$
Therefore we need to study the following sums:
\begin{prop}
Let $S_n=\sum_{i=1}^{n}\frac{\log p_i }{\sqrt{p_i}-1}$, $S'_n=\sum_{i=1}^{n}\frac{\log p_i }{p_i-1}.$
Then
\begin{align*}
 S_n\sim &2\sqrt{n\log{n}}\\
 S'_n\sim &\log{(n\log{n})}.
 \end{align*}
\end{prop}

\begin{preuve}
Let us introduce the function $\Lambda:\mathbb{N}\to\mathbb{R}$ defined by
$$\Lambda(n)=\begin{cases}
\ \log{p} & \text{ si } n=p \text{ is a prime number,}\\
\ 0 & \text{ otherwise.} \end{cases}$$

As $\vartheta(x)=\sum_{n\leq x}\Lambda(n)$, we have
$$\sum_{n\leq x}\frac{\Lambda(n)}{\sqrt{n}}=\frac{\vartheta{x}}{\sqrt{x}}+\int_2^x \frac{\vartheta(t)dt}{2t\sqrt{t}}.$$
Because $\vartheta(x)=x+o(x)$, the first term in the sum is equal to $\sqrt{x}+o(\sqrt{x})$ as $x\to\infty$.  Concerning the second, it is a divergent integral. As $\frac{\vartheta(x)}{2x\sqrt{x}}\sim \frac{1}{2\sqrt{x}}$, we find: $$\int_2^x \frac{\vartheta(t)dt}{2t\sqrt{t}}\sim\int_{2}^x\frac{dt}{2\sqrt{t}}\sim \sqrt{x}$$

As $p_n\sim n\log{n}$ and as $S_n$ is divergent, we get $$S_n\sim \sum_{i\leq n}\frac{\Lambda(i)}{\sqrt{i}}\sim 2\sqrt{n\log{n}}.$$

The same computation, using 
$$\sum_{p\leq x}\frac{\log p }{p-1}\sim \log x, \text{ as } x\to\infty,$$ leads to the estimate of $S'_n,$ \end{preuve}
\begin{coro}\label{ninv} For $i=1,2,$ 
$\delta_n^{(i)}\leq 1-\varepsilon^{(i)}$, where
\begin{align*}
\varepsilon^{(1)}\sim & \frac{8}{3\sqrt{n\log{n}}} \text{ and,}\\
\varepsilon^{(2)}\sim & \frac{4}{3n}.
\end{align*}
\end{coro}
\begin{preuve} 
We will prove it in the case of $\delta^{(1)}$ the second case being the same. As $\sum_{q}\phi_q\frac{\log q }{\sqrt{q}-1}\geq 2\frac{1}{g_{K_n}}S_n$, we get $$\delta_n\leq 1-2\frac{1}{g_{K_n}}S_n-(\gamma +\log{8\pi})\phi_\mathbb{C}(\mathcal{K_n})-(\gamma/2+\pi/4 + \log{\sqrt{8\pi}})\phi_\mathbb{R}(\mathcal{K}_n).$$
Then the result come from our estimates of the genus and of $S_n,$ and from the fact that archimedean factor are negligible. 
\end{preuve}

\section{On Prime Decomposition in Infinite Global Fields}

The class field theory produces examples of infinite global fields with finitely many places split. But as we can see from the last paragraph, it seems impossible to use it directly to obtain asymptotically good infinite global fields having infinitely many places split, and thus infinitely many positive invariants. The question of prime decomposition in infinite global field is central, and the analytic theory (in particular the Cebotarev density theorem) allows us to understand a bit what happens in such field. Let us recall a corollary of this theorem: 

\begin{prop}\label{pouet}
Let $K$ be a global field and let $T$ be the set of places of $K$ which split in a separable extension $L/K$ of degree $n$  such that, if $K$ is a function field, $\bar{\mathbb{F}_r}\cap L=\mathbb{F}_r$. Then $\bar{\delta}(T)\leq\frac{1}{n}$.
\end{prop}

\begin{coro}\label{ccebo} Let $K$ be a global field, $\KK/K$ an infinite global field  and $T$ the set of places of $K$ split in $\KK$. Then $\delta(T)$ exists and is equal to $0$.
\end{coro}

\subsection{Proof of Proposition \ref{pouet}}
Let us recall first the Cebotarev density theorem.
Let $K/L$ be a Galois extension, $G=Gal(L/K)$ and $\sigma\in G$. Put $$P_{L/K}(\sigma)=\left\{\mathfrak{p}\in P(K)\mid \ \ \ \exists\mathfrak{P}|\mathfrak{p}\in P(L) \ \ \ \sigma=\left(\frac{L/K}{\mathfrak{P}}\right)\right\},$$ and $<\sigma>=\left\{\tau^{-1}\sigma\tau \mid \  \tau\in G\right\}.$

\begin{theo} (Cebotarev)
The set $P_{L/K}$ has a Dirichlet density. It is given by $$\delta(P_{L/K}(\sigma))=\frac{\#<\sigma>}{\#G}.$$
\end{theo}
\begin{preuve} \cite[7.13.4]{NEU} 
\end{preuve}

Recall the following lemma:
\begin{lemme} Let L/K be a separable extension of global fields and $N$ be the Galois closure of $L/K.$ Then $\mathfrak{p}$ splits in $L/K$ if and only if it splits in $N/K.$
\end{lemme}
\begin{preuvelemma} $N/K$ is the compositum of all the conjugates of $L/K.$ If a place $\pg$ splits in $L/K,$ it also splits in $\sigma(L)/K$ for any $\sigma\in Gal(N/K),$ therefore it splits in their compositum. The converse is obvious.
\end{preuvelemma}

For any separable extension of global fields $L/K,$ define $P_{L/K}$ as the set of unramified places $\pg$ of $K$ admitting a place $\Pg$ of $L$ above it, whose inertia degree over $\pg$ is $1.$ 

\begin{lemme}Using notation and assumptions of the last lemma, put $G=Gal(N/K)$ and $H=Gal(N/L)$, One has : $$P_{L/K}=\bigsqcup_{<\sigma>\cap H\neq \emptyset}P_{N/K}(\sigma).$$
\end{lemme}
\begin{preuvelemma} \cite[7.13.5]{NEU}
\end{preuvelemma}

\begin{lemme}
Under this hypothesis, we have $\delta(P_{L/K})\geq \frac{1}{n}$ with the equality if and only if $L/K$ is Galois. In addition, we have $\displaystyle \delta(P_{L/K})\leq 1-\frac{n-1}{\#G}.$
\end{lemme}
In particular the set of places of $K$ having no degree $1$ places above them in $L$ is of positive Dirichlet density. Therefore we have the following result:
\begin{coro}\label{kunyavski} Let $\KK$ be an infinite number field (resp. function field). Let $U$ be the set of primes of $\Q$ (resp. $\FF$) such that $\phi_{\pg,\mathrm{N}\pg}=0.$ Then $U$ contains a set of positive Dirichlet density. 
\end{coro}

\begin{preuvelemma} For the first point see \cite[7.13.6]{NEU}. Concerning the second inequality, we start from: $$\delta(P_{L/K})=\delta\left(\bigsqcup_{<\sigma>\cap H\neq \emptyset}P_{N/K}(\sigma)\right)\leq \sum_{<\sigma>\cap H\neq \emptyset}\frac{\#<\sigma>}{\#G}.$$ Saying that $<\sigma>\cap H$ is not empty is the same as saying that $\sigma$ belongs to the set $$<H>=\bigcup_{\tau\in G} \tau H\tau^{-1},$$ and in this case $<\sigma>$ lies in $<H>.$ Therefore
$$\delta(P_{L/K})\leq \frac{1}{\#G}\#<H>.$$ An easy basic exercise in group theory (making $G$ act on $<H>$ by conjugation) shows that $$\#<H>\leq 1+\#G-[G:H],$$ from which we deduce the second inequality.
\end{preuvelemma}

\textit{Proof of Proposition \ref{pouet}:}  Consider the Galois closure $N/K$ of $L/K,$ and let $T$ be the set of splitting places of $L/K,$, or equivalently of $N/K.$ The last lemma implies that
 $$\bar{\delta}(T)\leq \delta(P_{N/K})=\frac{1}{[N:K]}\leq \frac{1}{n}.$$
 
To obtain the corollary, we just let $n\to +\infty$. \hfill $\square$

\subsection{Proof of Proposition \ref{iharaz}} 
We will prove in fact the following result:
\begin{prop} 
Let $\mathcal{K}$ be an infinite global field and $\{K_i\}_{i\in\mathbb{N}}$ be a tower representing $\KK.$ Let $T$ be the set of primes of $K_0$ split in $\KK/K_0.$ Then 
\begin{align*}
(NF) \quad \sum_{\mathfrak{p}\in T}\frac{\log{\mathrm{N}\mathfrak{p}}}{\mathrm{N}\mathfrak{p}-1}\leq & \frac{n_{K_0}^2}{\phi_\infty}(1-\frac{1}{2}(\log{2\pi}+\gamma)\phi_\infty).\\
(NF-GRH) \quad \sum_{\mathfrak{p}\in T}\frac{\log{{\mathrm{N}\mathfrak{p}}}}{\sqrt{\mathrm{N}\mathfrak{p}}-1}\leq & \frac{n_{K_0}^2}{\phi_\infty}(1-\frac{1}{2}\left(\log 8\pi+\gamma)\phi_\infty\right)\\
(FF) \quad \sum_{\mathfrak{p}\in T}\frac{\log_r{{\mathrm{N}\mathfrak{p}}}}{\sqrt{\mathrm{N}\mathfrak{p}}-1}\leq & \frac{n_{K_0}^2}{\phi_\infty}, 
\end{align*}
where the right hand terms can be infinite (in the case of asymptotically bad fields).
\end{prop}

\begin{preuve}
Let us make the proof in the number fields case, the function fields one and the GRH one being deduced using the appropriated inequalities. Suppose that $\KK$ is asymptotically good. 
Because of T-V Basic Inequalities, one has $$\sum_q \frac{\phi_q\log{q}}{q -1}+(\gamma/2+\log{2\sqrt{\pi}})\phi_\mathbb{R}+(\gamma+\log{2\pi})\phi_\mathbb{C}\leq 1.$$ 
As $\phi_\mathbb{R}+2\phi_\mathbb{C}=\phi_\infty,$ we obtain the inequality 
$$\sum_q \frac{\phi_q\log{q}}{q -1}\leq 1-\frac{1}{2}(\log{2\pi}+\gamma)\phi_\infty.$$ 

Let $T$ be the set $Dec(\mathcal{K}/K_0)$ of places of $K_0$ split in $\KK.$
For any $q$ in $A,$ let $T_q$ be the set of places of $T$ of norm $q,$ and let $\# T_q=N_q(T).$ Remark that, because of the Riemann-Hurwitz formula, the ratio $n_{K_i}/g_{K_i}$ is decreasing to $\phi_\infty.$ Therefore, for $q$ such that $T_q$ is not empty, 
$$\frac{\Phi_q(K_n)}{g(K_n)}\geq\frac{[K_n:K_0]}{g_{K_n}}=\frac{1}{n_{K_0}}\frac{n_{K_n}}{g_{K_n}}\geq \frac{1}{n_{K_0}}\phi_\infty.$$ Therefore $$\phi_q\geq \frac{1}{n_{K_0}}\phi_\infty.$$

Then $$\sum_{\mathfrak{p}\in T}\frac{\log{\mathrm{N}\mathfrak{p}}}{\mathrm{N}\mathfrak{p}-1}=\sum_q N_q(T)\frac{\log{q}}{q-1}\leq \sum_{T_q\neq\emptyset} n_{K_0}\frac{\log{q}}{q-1}\leq \frac{n_{K_0}^2}{\phi_\infty}\sum_{T_q\neq\emptyset} \phi_q\frac{\log{q}}{q-1},$$
$$\text{Thus }\ \ \ \sum_{\mathfrak{p}\in T}\frac{\log{\mathrm{N}\mathfrak{p}}}{\mathrm{N}\mathfrak{p}-1}\leq\frac{n_{K_0}^2}{\phi_\infty}\sum_{q} \phi_q\frac{\log{q}}{q-1}\leq  \frac{n_{K_0}^2}{\phi_\infty}(1-\frac{1}{2}(\log{2\pi}+\gamma)\phi_\infty).$$ 
\end{preuve}

\subsection{Asymptotically bad Case} One question arises naturally. We know that the Dirichlet density of the set of the split places is zero. When the field is asymptotically good, we have also proved that $s(T)$ has to be bounded. Can $s(T)$ be infinite in asymptotically bad towers? Can the set $T$ be infinite?  A negative answer would be a disaster, because of the link between the positive $\phi_q$ and the split places. We do not know at the moment how to construct asymptotically good fields having infinitely many positive $\phi_q$ or splitting places, and this seems to be a very difficult problem because of the restricted ramification condition. However, we can construct such examples in the asymptotically bad case: 

\begin{theo} \label{ex1} For any global field $K,$ there is an infinite Galois global field $\mathcal{K}$ containing $K$ such that infinitely many places split in $\mathcal{K}$, $\sum_{\mathfrak{p}\in Dec(\mathcal{K})}\mathrm{N} \mathfrak{p}^{-1}$ is infinite. It has no wild ramification over the Galois closure of $K$ (in some given separable closure of $K$) and $\delta(Ram(\mathcal{K}))$ is zero.
\end{theo}

 In the number field case (and also likely in the function field case) we could ask in addition that $s(Ram(\KK))\leq\epsilon$ (meaning: for any given $\epsilon>0,$ there is an infinite number field such that...) using a result of Gras (see \cite[Corollary $V.2.4.7$] {GRAS}) and doing the same as we will do in our proof. But he only proved it for the number fields case, that is why we will use once again Grunwald--Wang Theorem, which allows us to prove the theorem in both cases of number and function fields. 
 
\begin{preuve}
Recall that, for a set $E$ of places of a global field $K,$ we put $s(E):=\sum_{\mathfrak{p}\in E}\mathrm{N} \mathfrak{p}^{-1},$ eventually infinite.
Let $K_0$ be the Galois closure of $K.$ Consider, for the simplicity of notation, the number field case. Let $S_0$ (resp. $D_0$) be the ramification locus of $K_0/\mathbb{Q}$ (resp. the splitting locus of $K_0/\mathbb{Q}$). We will extend this notation to other indices than $0$ by putting $S_n=Ram(K_n/\Q)$ and $D_n=Dec(K_n/\mathbb{Q})$. Let $T_0$ be a finite set of places of $\mathbb{Q}$ split in $K_0/\mathbb{Q}$ such that $s(T_0)\geq 1.$ Such a set exists since $\delta(D_0)>0$ by the Cebotarev density theorem. Let $\ell$ be a prime in $T_0$ (or a prime number different from the characteristic of  $K_0$ in the function field case).

We construct by induction the tower of fields $K_n,$ Galois over $\Q$ and the set $T_n\subset D_n$ with $s(T_n)\geq n+1$ having the desired properties in the following way. Suppose that we have constructed a field $K_{n-1},$ $n\geq 1,$ Galois over $Q,$ and a set $T_{n-1}$ as above. Let $K_n^\ast$ be a $\ell$-extension of $K_{n-1}$ unramfied outside of $D_{n-1}(K_{n-1}),$ such that all the places above $T_{n-1}$ in $K_{n-1}$ split. Let $S_n^\ast$ denotes the ramified places of $K_n^\ast/K_{n-1}.$ Then we consider the maximal $\ell$-extension $K_n$ of $K_{n-1},$ unramified outside of the places of $S_n^\ast$ and their conjugates by the Galois action, where the places of $T_{n-1}(K_{n-1})$ split in $K_n.$ This extension is non trivial, and moreover $K_n$ is Galois over $\Q.$ To see that, let us take a morphism $\sigma$ from $K_n$ to a separable closure $\bar{\Q}$ of $\Q.$ $\sigma(K_n)$ is a $\ell$-extension of $\sigma(K_{n-1})=K_{n-1}$ unramified outside of the conjugates of places of $S_n^\ast,$ where the places above $T_{n-1}$ split, therefore it is contained in $K_n.$ For the set $T_n,$ let us take a subset of $D_n$ containing $T_{n-1}$ and such that $s(T_n)\geq n+1.$

If we denote by $T$ the splitting locus of $\mathcal{K}/\mathbb{Q}$ and put $S=\cup S_n$ for the ramification locus of $\mathcal{K}/\mathbb{Q},$ then we have $s(T)\geq s(T_n)$ for any $n$ and therefore $s(T)=\infty.$ In addition, $S_{n+1}$ is contained in $D_{n}$ for any $n\geq 0.$ Therefore $\cup_{m>n}S_m\subset D_n$ for any $n\geq 0.$ So we have $\delta(S)=0.$ Indeed, suppose that $\bar{\delta}(S)>1/\ell^{n_0}$ for a given $n_0>0.$ Then $\bar{\delta}(S-\cup_{i\leq n_0}S_i)>1/\ell^{n_0},$ but this is impossible because $S-\cup_{i\leq n_0}S_i\subset D_{n_0}$ and $\delta(D_{n_0})\leq1/\ell^{n_0}.$ $\ell$ is not ramified in $\mathcal{K}/\Q,$ thus there is no wild ramification over $K_0.$
\end{preuve}
 
To conclude, we remark that there are infinite global field having no places split: this is for example the case for the maximal abelian extension of $\Q.$

\subsection{Abelian Case} Let us give now the negative answer to a natural conjecture of Michel Balazard formulated at Poncelet Laboratory seminar (Moscow). In the case of $\mathbb{Q}$ and in that  of cyclotomic extensions, we have the following beautiful result of Norton \cite{Nor}, which can be obtained using Br\"un--Titchmarch and Siegel--Walfisz Theorems:

\begin{prop}[Norton] Let $q$ be a prime number and $a$ be an integer non divisible by $q.$ Let $$I_x=\{p \text{ prime number } \vert\, p=a\, \mod q, p\leq x \}.$$ Then there exists a constant $M$ independent from $q,a$ such that $$\bigg|\sum_{p\in I_x}\frac{1}{p}-\frac{1}{q-1}\log\log{x}\bigg|\leq M.$$
\end{prop}

It seems to be hard to obtain good generalisations of these two theorems to general global fields (see \cite{GSW} for the number fields case).
However Michel Balazard put forward the following conjecture:
\begin{conjecture} \label{bal} Let $K$ be a global field. Let $L/K$ be an abelian extension of degree $n.$ Then there exists a constant $M$ depending only on $K$ such that
$$\bigg|\sum_{\mathfrak{p}\in Dec({L/K}),\ \mathrm{N} \mathfrak{p}\leq x}\frac{1}{\mathrm{N} \mathfrak{p}}-\frac{1}{n}\log\log{x}\bigg|\leq M.$$
\end{conjecture}
Unfortunately, this conjecture is not true, at least this formulation, as we will see it producing an example of pro-cyclic infinite global field which violates it. 

\begin{prop}\label{bala} Both in the case of number fields and in the case of function fields, there exists an infinite Galois  global field $\KK/K$ of pro-cyclic Galois group such that  $\sum_{\mathfrak{p}\in Dec(\mathcal{K})}\frac{1}{\mathrm{N} \mathfrak{p}}=\infty.$
\end{prop}
\begin{preuve}
Grunwald-Wang Theorem will allow us to produce such an example. We will give the construction in the number field case for simplicity of notation. Let us take a prime number $p_0=3$ and consider a cyclic extension $K_0$ of the field $K$ ($K$ being a finite separable extension of $\mathbb{Q})$ of degree $p_0.$ 

Let us take then a finite set $T_0$ of primes $\mathfrak{p}_1^0,...,\mathfrak{p}_n^0$ of the ground field $\Q$ split in $K_0$ such that $\sum_{\mathfrak{p}\in T_0}\frac{1}{\mathrm{N} \mathfrak{p}}\geq 1.$ 

Take for the field $K_1$ a cyclic extension of $\Q$ of prime degree $p_1>p_0$ which does not divide the degree of $K,$ such that $T_0$ is split in $K_1.$ Thus $K_1.K_0/K$ is a cyclic extension of degree  $p_1p_0,$ where $T_0$ is split. Take for $T_1$ a set of places of $\Q$ split in $K_1.K_0,$ containing $T_0$ and satisfying $\sum_{\mathfrak{p}\in T_1}\frac{1}{\mathrm{N} \mathfrak{p}}\geq 1.$

Suppose that we have constructed in this way a field $K_n.K_{n-1}...K_0/K$ of degree $p_n...p_0$ together with a set $T_n$ of splitting places, satisfying $\sum_{\mathfrak{p}\in T_n}\frac{1}{\mathrm{N} \mathfrak{p}}\geq n.$ We construct $K_{n+1},T_{n+1}$ from $K_n.K_{n-1}...K_0,T_n$ in the same way that we have constructed $K_1,T_1$ from $K_0,T_0.$ Then we obtain our tower by induction. This tower satisfies the hypotheses of the proposition.
\end{preuve}

\begin{coro}\label{cbal} Conjcture \ref{bal} is false, even for cyclic extensions.
\end{coro}
\begin{preuve} Suppose it is true for a given $K$. Consider the tower of the last proposition starting from $K$, and let  $T$ be its splitting locus. Taking the limit as $n\to+\infty,$ for any $x$ we have: $$\sum_{\mathfrak{p}\in T,\ \mathrm{N} \mathfrak{p}\leq x}\frac{1}{\mathrm{N} \mathfrak{p}}\leq M.$$ Taking the limit as $x\to+\infty$ we get a contradiction.
\end{preuve}

\section{Proof of Theorem \ref{phi2}}\label{schmi}

First let us recall the following definitions and notation from Galois cohomology. Let G be a profinite group. The cohomological dimension $cd\, G$ of $G$ is the smallest integer $n,$ if it exists, such that $$H^q(G,A)=0 \text{ for every } q>n $$ and every torsion $G$-module $A.$ If no such integer $n$ exists, then $cd\,G=\infty.$ For a pro-$p$-group, $cd\,G\leq n$ if and only if $H^{n+1}(G,\mathbb{Z}/p\mathbb{Z})=0$ (see \cite[$\S$III.$3$]{NCG}). For example, any non trivial finite group has an infinite cohomological dimension.

We will prove now the following result which implies Theorem \ref{phi2}.
\begin{theo}\label{invnul} Let $p$ be an odd prime number and let $S_0$ be a set of primes congruent to $1$ modulo $p,$ such that $\Q_{S_0}(p)$ --- the maximal $p$-extension of $\Q$ unramified outside of $S_0$--- is of cohomological dimension $2.$ Let $P$ be a finite set of prime numbers. Then there exists a finite set of primes $S$ containing $S_0$ and not containing $p$ such that, for any $m>1$ and any prime number $\ell,$ $\phi_{\ell^m}(Q_S(p))=0$, and that, for any $\ell\in S\cup P,$ $\phi_{\ell}(\Q_S(p))=0.$
\end{theo}
Remark that the only primes which can be ramified in a $p$-extension are congruent to $0,1$ modulo $p,$ thus the condition on $S_0$ is natural.

 As a corollary, we deduce directly Theorem \ref{phi2}. Indeed, Labute (see \cite{LMG}) gave some examples having the cohomological dimension $2$ property, for example $cd\, Gal(\Q_{\{7,19,61,163\}}(3)/\Q)=2,$ so one can apply Theorem \ref{invnul} to the set $P.$ The resulting infinite number field $\KK=\Q_S(3)$ is asymptotically good because it is tamely ramfied and unramified outside of a finite set of primes. In addition, the splitting prime numbers $\ell$ in $\KK$ satisfy $\phi_{\ell}>0$ because of Proposition \ref{propdec}, and no other prime number belongs to $\mathrm{PSupp}(\KK/\Q)$ because of Theorem \ref{invnul}.

Let us now prove Theorem \ref{invnul}. In order to do that, recall a result due to Schmidt (\cite[Theorems $2.2$ and $2.3$]{SCS}). Denote by $G_S(p)$ the Galois group of the maximal  $p$-extension of $\Q$ unramified outside of $S.$
\begin{theo} [Schmidt]\label{schmidt} Let $p$ be an odd prime number and let $S$ be a finite set of prime numbers congruent to $1$ modulo $p.$
\begin{enumerate}\item Suppose that $G_S(p)\neq 1$ and that $cd\, G_S(p)\leq 2.$ Then $cd\, G_S(p)=2,$ and for any $\ell\in S,$ $\Q_S(p)$ realises the maximal $p$-extension of $\Q_\ell.$ 
\item Suppose $cd\, G_S(p)=2.$ Then, if $\ell\notin S$ is an other prime number congruent to $1$ modulo $p,$ which is not split in $\Q_S(p)/\Q.$ Then $cd\, G_{S\cup\{\ell\}}(p)=2.$ \end{enumerate}
\end{theo}

Consider now an odd prime number $p,$ $S_0$ a set of primes congruent to $1$ modulo $p,$ such that $cd\, G_{S_0}=2$ and a finite set $P=\{p_1,...,p_r\}$ of prime numbers.
Let us begin by proving two lemmas:

\begin{lemme}
If $cd\, G_S(p)=2,$then for any $m>1,$ for any prime number $\ell$ not belonging to $S,$ $\phi_{\ell^ m}(\Q_S(p))=0.$
\end{lemme}
\begin{preuvelemma}  If $\mathrm{cd}G_S(p)=2,$ then for any closed subgroup $H$ of $G_S(p),$ $H^n(H)=0$ for any $n>2$ (\cite[$\S$I prop. $21'$]{SCG}). Suppose that $\phi_{\ell^m}>0$ for $\ell\notin S$ and $m>1.$ Then the Frobenius of $\mathbb{F}_{\ell^m}$ which is of finite order, can be lifted in the decomposition subgroup in $G_S(p)$ of any place over $\ell$  to an element of finite order. The subgroup generated by it is finite therefore of infinite cohomological dimension, so it does not satisfy the above property. This leads to a contradiction.
\end{preuvelemma}

\begin{lemme}If $cd\, G_S(p)=2,$ then for any $m,$ for any prime $\ell$ in $S$, $\phi_{\ell^m}(\Q_S(p))=0.$
\end{lemme}
\begin{preuvelemma} We apply Theorem \ref{schmidt}, which tells us that the inertia index of $\ell$ in $\Q_S(p)$ is infinite (indeed, there exists an unramified infinite algebraic extension of $\Q_\ell$).
\end{preuvelemma}

Let us prove now Theorem \ref{invnul}.

\begin{preuve} Consider $\Q_{S_0}(p)$ and its splitting locus $T.$ We will complete the set $S_0$ by a finite number of places in order to make the $\phi_p$ vanish one after the other. We will construct a set $S$ of primes having the following properties:

\begin{enumerate}
\item any $q\in S$ is congruent to $1$ modulo $p,$
\item $\Q_{S}(p)$ is asymptotically good and $cd\,G_{S}(p)=2,$
\item $p_i,$ $i\leq r,$ is not split in $\Q_{S}(p).$
\end{enumerate}

The key ingredient of the proof is the existence of a tamely ramified $p$-extension, unramified at $T,$ such that the $p_i$ are inert. Such an extension exists because of the precise version of  Grunwald--Wang Theorem.

Indeed, because of the Grunwald--Wang theorem, there is a cyclic extension $K$ of degree $p$, unramfied outside of $(P(\Q)-T),$ such that $p$  and the $p_i$ are inert. Put $S=S_0\cup\mathrm{Ram}(K).$ Thus no $p_i$ is split in $\Q_S(p).$

As the only primes which can be ramified in a $p$-extension are congruent to $0,1$ modulo $p,$ $S$ satisfy the first condition. Moreover $\Q_{S}(p)$ is asymptotically good, because it is tamely and finitely ramified.
In order to prove ii, we have to prove that one can apply Theorem \ref{schmidt} (ii) to all the places in  $\mathrm{Ram}(K)=\{s_1,...,s_m\}.$ We have to satisfy that $s_i$ is not split in $\Q_{S_0\cup \{s_0\}\dots\cup \{s_{i-1}\}}(p).$ As this field contains $\Q_{S_0},$ and as $s_i$ is not split in it since $T$ is unramified, $s_i$ is not split in $\Q_{S_0\cup \{s_0\}\dots\cup \{s_{i-1}\}}(p).$ Thus we can apply Theorem \ref{schmidt} (ii) to all the set $\mathrm{Ram}(K)$ and obtain (ii).

From the first lemma, we deduce that the $\phi_{p_i^m}(\Q_{S})$ are all zero and we have proven the theorem.   \end{preuve}

Let us conclude this section by the following remark: in order to prove Theorem \ref{phi2} we could have first considered a $\ell$-extension $K/\Q$ where all the prime in $P$ are inert and its ramification locus $\mathrm{Ram}(K)$  and then extended this set to a set $S$ such that $cd\,G_S(p)=2$ (see \cite[7.3]{SCS}). Using this, the ramification locus and thus its deficiency may be much bigger.


\bibliographystyle{amsplain}
\bibliography{biblio}

\end{document}